\documentclass[onecolumn, 10 pt, conference]{cssconf}  

\IEEEoverridecommandlockouts                              
\newtheorem{theorem}{Theorem}
\newtheorem{corollary}{Corollary}
\newtheorem{definition}{Definition}

\newtheorem{proposition}{Proposition}
\newtheorem{assumption}{Assumption}
\newtheorem{lemma}{Lemma}
\newtheorem{remark}{Remark}

\usepackage{amsmath} 
\usepackage{amssymb}  
\usepackage{mathtools}
\usepackage{mathrsfs}
\usepackage{braket}
\usepackage{pgfplots}
\usepackage{subcaption}
\usepackage[utf8]{inputenc}
\usepackage[backend=biber,style=ieee, url=false, isbn=false,doi=false]{biblatex}
\DeclareFieldFormat{doi}{%
  \hfil\penalty99\hfilneg\space {\textsc{\MakeLowercase{DOI}}} \addcolon\addnbspace
  \ifhyperref
    {\href{https://doi.org/#1}{\nolinkurl{#1}}}
    {\nolinkurl{#1}}}
\pgfplotsset{compat=1.17}

\usepackage{tikz}
\usepackage{xcolor}

\newcommand{\J}{\mathscr{J}}

\usepackage[hidelinks,colorlinks,linkcolor=blue,citecolor=green!70!black]{hyperref}
\usepackage{xurl}
\renewcommand{\theenumi}{\alph{enumi}}

\pgfplotsset{
scriptsize/.style={
width=4.5cm,
height=,
legend style={font=\scriptsize},
tick label style={font=\scriptsize},
label style={font=\footnotesize},
title style={font=\footnotesize},
every axis title shift=0pt,
max space between ticks=13,
every mark/.append style={mark size=8},
major tick length=0.1cm,
minor tick length=0.066cm,
},
}
\usetikzlibrary{intersections}
\usepgfplotslibrary{groupplots}

\pgfplotsset{
small/.style={
width=6.3cm,
height=,
tick label style={font=\footnotesize},
label style={font=\small},
max space between ticks=25,
},
}

\title{\LARGE \bf
Probabilistic Reachable Set Estimation for Saturated Systems with Unbounded Additive Disturbances
}

\author{Carlo Karam, Matteo Tacchi-Bénard, Mirko Fiacchini%
\thanks{C. Karam, M. Tacchi-Bénard and M. Fiacchini are with Univ. Grenoble Alpes, CNRS, Grenoble INP (Institute of Engineering
Univ. Grenoble Alpes), GIPSA-lab, 38000, Grenoble, France.  
\{\tt carlo.karam, matteo.tacchi, mirko.fiacchini\}\tt @gipsa-lab.fr.}
}

\begin{document}

\newcommand{\R}{\mathbb{R}}
\newcommand{\N}{\mathbb{N}}
\newcommand{\E}{\mathbb{E}}
\newcommand{\Tr}{\mathrm{Tr}}
\newcommand{\Ek}[1]{\E\left[e_{#1}^\top P e_{#1}\right]}
\newcommand{\Eq}[1]{\E\left[q_{#1}\right]}
\newcommand{\ABKc}{\left(A + \sum_{i \in \mathscr{J}} B_{(i)}K_{i}\right)}

\maketitle
\thispagestyle{empty}
\pagestyle{empty}

\begin{abstract}
\textbf{In this paper, we present an analytical approach for the synthesis of ellipsoidal probabilistic reachable sets of saturated systems subject to unbounded additive noise. Using convex optimization methods, we compute a contraction factor for the saturating error dynamics, which allows us to tightly bound their evolution and thereby construct accurate reachable sets. The proposed approach is applicable to independent, zero mean disturbances with a known covariance. A numerical example illustrates the applicability and effectiveness of the proposed design.}
\end{abstract}

\section{Introduction}
Reachability analysis is of particular relevance for the increasingly popular field of
stochastic model predictive control (SMPC) \cite{mesbahStochastic2016}, as it enables
the characterization of the reachable sets of a dynamical system subject to stochastic
disturbances. This is essential for providing probabilistic guarantees on key properties
such as constraint satisfaction, recursive feasibility, and closed-loop stability.

The problem of constructing probabilistic reachable sets (PRS) has been extensively studied in the case of bounded noise, with well-established methods available in the literature, e.g. \cite{bonzaniniTubebased2019,schluterConstraintTightening2020,sasfiRobust2023,cannonRobust2010}. However, the
synthesis of such sets under unbounded probabilistic noise remains an open problem. Data-driven 
and sampling-based approaches have proven effective in constructing PRS under unbounded 
disturbances, as in \cite{kerzDataDriven2023, panDatadriven, hewingScenarioBased2020, 
fochesatoDataDriven2022}, but these methods rely on access to noise samples. In contrast, 
purely analytical approaches have been proposed \cite{hewingStochastic2018,hewingRecursively2020}, 
with recent advances incorporating more realistic assumptions, such as sub-gaussian 
\cite{aoStochastic2025} and time-correlated \cite{fiacchiniProbabilistic2021} disturbances. 
A key limitation of these approaches is their treatment of chance constraints on both states 
and inputs, whereas real-world applications often require hard input constraints to 
properly account for actuator saturation. Some existing results address hard input constraints 
by relaxing guaranteed feasibility and parameterizing the control input as a function of 
past disturbances \cite{hokayemStochastic2009, paulsonStochastic2020}. However, these methods 
face computational challenges in high-dimensional settings and may be affected by relevant conservatism. For output-feedback with hard input constraints, 
\cite{joaOutput2023} derives uniform time-invariant bounds under Gaussian noise assumptions, 
successfully constructing conservative PRS. Very recently, \cite{kohlerPredictive2025} proposed 
a contraction-based approach for constructing PRS in the presence of hard input constraints, 
ensuring closed-loop SMPC properties through open-loop system analysis.

This work presents an approach for synthesizing probabilistic reachable sets in the context
of an indirect feedback formulation, where the nominal and error dynamics are decoupled. We consider linear systems 
affected by unbounded, independent, zero-mean noise, and subject to hard input constraints modeled 
as a direct saturation on the control input. Leveraging established results from saturated systems theory \cite{tarbouriechStability2011,fiacchiniQuadratic2012}, we bound the saturated error dynamics within a well-defined convex set. This set is constructed such that its vertices encode all possible saturation scenarios: fully saturated (open-loop), unsaturated (closed-loop), and partially saturated configurations. This bounding enables the computation of quadratic Lyapunov functions defined on that set, which remain valid for the error dynamics. Motivated by the realization that the system probabilistically switches between
saturation modes, we compute its effective contraction rate which lies in between the open-loop
(worst-case) and unsaturated closed-loop (best-case) rates. We subsequently derive a tight
bound on the expectation of a quadratic transformation on the error, enabling the computation
of accurate ellipsoidal probabilistic reachable sets with user-defined violation probability. To the best of our knowledge, no constructive PRS design method incorporating the feedback on the error term exists in the nonlinear (saturated) setting. 

This novel PRS design technique
is particularly effective at reducing the conservatism of the existing approaches without 
requiring extensive assumptions on the noise. The results herein have significant implications 
for reachability-based SMPC frameworks, as the obtained sets can be used for an appropriate constraint 
tightening that ensures probabilistic recursive feasibility. For this reason, we frame our development within a
stochastic predictive control context, but leave the actual SMPC implementation for a future work.

\subsection{Notation}
Throughout the article, $\mathbb{R}$ denotes the set of reals, $\mathbb{N}_0$ denotes the set of nonnegative integers, and $\mathbb{N}^+_m$ denotes the set of positive integers from 1 to $m$. Variables at time step $k$ are indexed with the corresponding subscript, e.g. $x_k$. The columns of a matrix $A$ are denoted as $A_{(i)}$, while its rows are denoted $A_{i}$. Given $P \succ 0$, define the ellipsoid $\mathscr{E}(P,r) = \{x \in \R^n | x^\top P x \leq r\}$ where $r \geq 0$ is the ellipsoidal scaling parameter. A random variable $e$ following a distribution $Q^{e}$ is denoted by $e \sim \mathcal{Q}^e$. The powerset $\mathscr{J}$ of a set $\mathcal{S}$ is denoted as $\mathcal{J} \subseteq \mathcal{S}$. Finally, the symbol $\left\| \cdot \right\|_\infty $ denotes the $\ell^\infty$-norm of a vector, i.e., $\left\| \cdot \right\|_\infty = \max_{i} \, \left| \cdot \right|_{i} $.

\section{Preliminaries}
In this section, we present the theoretical framework which constitutes the foundation of our main development.

\subsection{Definitions}\label{sec:def}
We introduce key concepts related to probabilistic reachable sets that will be central in our
characterization of the system (\ref{eq:lti-sat-system}).

\begin{definition}[Probabilistic Reachable Set] 
   A sequence of sets $\mathcal{R}_k$ are said to be probabilistic reachable sets (PRS) of violation level
   $\varepsilon \in [0, 1]$ for a stochastic process $\left\{e_k\right\}_{k \in \mathbb{N}_0}$ 
   taking values in $\mathbb{R}^n$ if
   \begin{equation*}
      e_0 = 0 \implies \Pr\left\{ e_k \in \mathcal{R}_k \right\} \geq 1 - \varepsilon, 
   \quad \forall k > 0 \end{equation*}
\end{definition}

\begin{definition}[Probabilistic Ultimate Bound, \cite{kofmanProbabilistic2012}]
   A set $\mathcal{R}$ is said to be a probabilistic ultimate bound (PUB) of violation level
   $\varepsilon \in [0, 1]$ for a stochastic process $\left\{e_k\right\}_{k \in \mathbb{N}_0}$ 
   taking values in $\mathbb{R}^n$, if
   \begin{equation*}
      \forall e_0 \in \mathbb{R}^n, \exists t \,\, \text{s.t.}\,\, \Pr\left\{e_k \in \mathcal{R} \right\} \geq 1 - \varepsilon 
      \quad \forall k \geq t.
   \end{equation*}
\end{definition}

\begin{definition}[Probabilistic Invariant Set]
   A set $\mathcal{R}_\infty$ is said to be a probabilistic invariant set (PIS) of violation level
   $\varepsilon \in [0, 1]$ for a stochastic process $\left\{e_k\right\}_{k \in \mathbb{N}_0}$ 
   taking values in $\mathbb{R}^n$, if
   \begin{equation*}
      e_0 \in \mathcal{R}_\infty \implies \Pr\left\{e_k \in \mathcal{R}_\infty \right\} \geq 1 - \varepsilon 
      \quad \forall k > 0.
   \end{equation*}
\end{definition}

\subsection{Problem Statement}
Consider the linear time-invariant (LTI) system subject to additive disturbances, whose 
discrete-time dynamics are given by
\begin{equation}\label{eq:lti-system}
   x_{k+1} = Ax_k + Bu_k + w_k,
\end{equation}
with state $x_k \in \mathbb{R}^{n}$, inputs $u_k \in \mathbb{R}^{m}$, process noise
$w_k \in \mathbb{R}^{n}$, and discrete time $k \in \mathbb{N}_0$. The matrices $A \in \mathbb{R}^{n \times n}$
and $B \in \mathbb{R}^{n \times m}$ are known, and the state $x_k$ is measurable. 

\begin{assumption}\label{ass:noise-iid}
The noise $w_k \sim \mathcal{Q}^w$ is independently and identically distributed (i.i.d.) with $\E\left[w_k\right] = 0$ and $\E\left[w_k w_k^\top \right] = W \succeq 0$, for all $k \in \mathbb{N}_0$. 
\end{assumption}
 
Notice that no knowledge of the distribution $\mathcal{Q}^w$ is required outside of its mean and variance. 
The control problem typically studied in the context of stochastic MPC enforces chance constraints on the state with violation level $\varepsilon_x$, i.e.,
\begin{equation}\label{eq:state-chance-constraints}
   \Pr\left\{ x_k \in \mathcal{X} \right\} \geq 1 - \varepsilon_x \quad \forall k \in \mathbb{N}_0,
\end{equation}
where $\mathcal{X}$ is a convex set, and hard input constraints 

\begin{equation}\label{eq:input-constraints}
   u_k \in \mathcal{U}\quad \forall k \in \mathbb{N}_0,
\end{equation}
where $\mathcal{U}$ is a convex set. To simplify the synthesis of reachability-based control laws, we model (\ref{eq:input-constraints}) as input saturations. This approach allows us to leverage established results from the theory of saturated systems to compute reachable sets while remaining faithful to conditions encountered in practical applications (i.e., actuator saturation). This formulation is applicable to any symmetric polytopic constraint centered at the origin. The considered system thus becomes
\begin{equation}\label{eq:lti-sat-system}
   x_{k+1} = Ax_k + B\varphi\left(u_k\right) + w_k,
\end{equation}
subject to (\ref{eq:state-chance-constraints}), where $\varphi : \mathbb{R}^m \to \mathbb{R}^m$ 
is the saturation function defined as
\begin{equation}\label{eq:sat-func}
   \varphi_{i}(u) = \mathrm{sign}\left(u_{i}\right) \min\left\{\left| u_{i} \right|, 1 \right\}
   \quad \forall i \in \mathbb{N}^+_m.
\end{equation}
Without loss of generality, the saturation bounds are assumed to be equal to 1. 

Reachability-based control approaches consist in splitting the system state $x_k$ 
into nominal and error parts
\begin{equation}
   x_k = z_k + e_k,
\end{equation}
such that the error $e_k$ is kept in a bounded neighborhood of the nominal trajectory $z_k$ 
via a locally-stabilizing feedback controller. As such, the input to (\ref{eq:lti-sat-system}) becomes
\begin{equation}\label{eq:tube-based-input}
   u_k = v_k + Ke_k,
\end{equation}
where $v_k$ is the control action computed by the nominal controller (i.e., MPC) designed for $z_k$. Defining
\begin{equation}
   f(e_k, v_k) = Ae_k + B \left(\varphi(Ke_k + v_k) -v_k\right), 
\end{equation}
we get the following discrete-time dynamics for the nominal trajectory and the error
\begin{subequations}\label{eq:split-dynamics}
   \begin{align}
      z_{k+1} &= Az_k + Bv_k,\label{eq:nominal-dynamics}\\
      e_{k+1} &= f(e_k, v_k) + w_k,\label{eq:error-dynamics}
   \end{align}
\end{subequations}
with $z_0 = x_0$ and $e_0 = 0$.

We can ensure that the uncertain system satisfies the state chance constraints (\ref{eq:state-chance-constraints})
for all disturbance realizations $w_k$, if the nominal trajectory dynamics (\ref{eq:nominal-dynamics}) are
appropriately bounded \cite{rawlingsModel2017}. These bounds can be obtained through a tightening of 
the original constraint set $\mathcal{X}$ by a set $\mathcal{R}_k$, the PRS of the uncertain system 
(\ref{eq:lti-sat-system}) (see \cite{hewingStochastic2018}). Given the error dynamics (\ref{eq:error-dynamics}),
our goal is to determine such a set $\mathcal{R}_k$, which we construct in what follows by estimating an 
ellipsoidal confidence region for $e_k$. A similar tightening is typically required on the input constraint set
$\mathcal{U}$. However, this becomes unnecessary when (\ref{eq:input-constraints}) are directly integrated into the
dynamics as saturations, which guarantees their satisfaction by design (see Remark \ref{rem:u-tightening}). 
For our development, the following regularity assumptions are considered.

\begin{assumption}\label{ass:regularity-conds}\leavevmode
   
   \begin{enumerate}
      \item The matrix $A$ is Schur.
      \item The nominal input $\left\| v_k \right\|_\infty \leq 1 $, $\forall k \in \mathbb{N}_0$. \label{cond:v-infty}
   \end{enumerate}
\end{assumption}

\begin{remark}\label{rem:u-tightening}
   Condition \ref{cond:v-infty} on the nominal control action is equivalent to imposing the input constraint set $\mathcal{U}$.
\end{remark}


\section{Probabilistic Reachable Set Synthesis}

In this section, we present our approach for synthesizing a sequence of PRS for the stochastic saturated system (\ref{eq:error-dynamics}).

\subsection{Convex Bounds of Saturated Functions}\label{sec:co-bounds-sat}
Following the same approach as in \cite{fiacchiniQuadratic2012}, we will first show that the image of the error $e_k$ through
the saturated function $f(e_k, v_k)$ is contained within a well-defined set. This will eventually allow us 
to obtain an upper-bound on the PRS of the uncertain system.

\begin{definition}[Support Function]
   The support function of a nonempty convex set $C \subseteq \mathbb{R}^n$, at $\eta \in \mathbb{R}^n$ is given by
   $\displaystyle \delta_C(\eta) = \sup_{x \in C} \, \eta^\top x$.
\end{definition}

The properties of support functions are such that the following set inclusion theorem holds.

\begin{theorem}[\cite{rockafellarConvex2015}]\label{th:support-inclusion}
   \emph{Given nonempty convex closed sets $C_1$ and $C_2$, $C_1 \subseteq C_2$ if and only if 
   $\delta_{C_1}(\eta) \leq \delta_{C_2}(\eta)$ for all $\eta \in \mathbb{R}^n$.}
\end{theorem}

Using the above result, we present the following theorem, which holds for all $e \in \mathbb{R}^n.$

\begin{theorem}\label{th:sat-incl-co}
   \emph{If $v$ is such that $\left\| v \right\|_\infty \leq 1$, then the function 
   $f(e, v) = Ae + B \left(\varphi(Ke + v) - v\right) \in F(e)$ with
   \begin{equation}\label{eq:co}
      F(e) = \mathrm{co}\left(\left\{  \Big(A + \sum_{i \in \mathscr{J}} B_{(i)}K_i \Big) e  \right\}_{\mathscr{J} \subseteq \mathbb{N}^+_m} \right).
   \end{equation}
} 
\end{theorem}

\proof First, assume $m = 1$ and suppose that $Ke \geq 0$. Since $-1 \leq v \leq 1$, we have 
$-1 \leq Ke + v \leq 1 + Ke$, hence
\begin{equation*}
   \varphi(Ke + v) = 
   \begin{cases}
      Ke + v & \text{if } Ke + v \leq 1 \\
      1 & \text{if } Ke + v \geq 1  \\
   \end{cases}
\end{equation*}
Subsequently,
\begin{equation*}
   \varphi(Ke + v) - v = 
   \begin{cases}
      Ke & \text{if } Ke + v \leq 1 \\
      1 - v& \text{if } Ke + v \geq 1  \\
   \end{cases}
\end{equation*}
Since $Ke$ is unbounded above and $1 - v \geq 0$, we have $0 \leq \varphi(Ke + v) - v \leq Ke$.
Analogously, we can show that if $Ke < 0$, we have $Ke \leq \varphi(Ke + v) - v \leq 0$. Therefore,
for arbitrary $m$, $\mathscr{J} \subseteq \mathbb{N}_{m}^+$, $i \in \mathscr{J}$, we have
\begin{equation*}
\begin{cases}
    \hphantom{K_ie}0 \leq \varphi \left( K_ie + v_i \right) - v_i \leq K_{i}e & \quad \text{if } Ke \geq 0, \\
    \hphantom{0}K_{i}e \leq \varphi \left( K_ie + v_i \right) - v_i \leq 0 & \quad \text{if } Ke < 0.
\end{cases}
\end{equation*}
This implies that
\begin{equation*}
   \eta^\top B_{(i)} \left( \varphi(K_ie + v_i) - v_i \right) \in 
   \mathrm{co}\left\{ 0, \eta^\top B_{(i)} K_{i} e \right\}
   \subseteq \mathbb{R},
\end{equation*}
holds for all $\eta \in \mathbb{R}^n$ and all $i \in \mathscr{J}$. Equivalently, we have
\begin{align*}
   \min\left\{ 0, \eta^\top B_{(i)} K_{i} e \right\} &\leq 
   \eta^\top B_{(i)} \left( \varphi(K_ie + v_i) - v_i \right)\\ 
   &\leq \max\left\{ 0, \eta^\top B_{(i)} K_{i} e \right\}.
\end{align*}
Thus, given $\eta \in \mathbb{R}^n$, $\|v\|_{\infty} \leq  1$, and $e \in \R^n$, there exists 
$\mathscr{J}(e, \eta) \subseteq \mathbb{N}^+_m$ such that
\begin{align*}
   \begin{split}
      \eta^\top f(e, v) &= \eta^\top Ae + \sum_{i \in \mathbb{N}^+_m} \eta^\top B_{(i)} \left(\varphi(K_ie + v_i) - v_i \right) \\
      &\leq \eta^\top Ae + \sum_{i \in \mathscr{J}(e, \eta)} \eta^\top B_{(i)} \left(\varphi(K_ie + v_i) - v_i \right),
   \end{split}
\end{align*}
holds. From Theorem \ref{th:support-inclusion}, we have $f(e, v) \in F(e)$.
\endproof

Theorem \ref{th:sat-incl-co} states that for all $e \in \R^n$ and any $\left\| v \right\|_\infty \leq 1$, the
image $f(e, v)$ is contained in a known polytope $F(e)$ whose vertices are given by the linear combination $\ABKc e$, where $\mathscr{J} \subseteq \mathbb{N}^+_m$. This means that the polytope vertices encode all possible saturation scenarios: full, partial, and no saturation.

\begin{remark}
   In this work, the inclusion in Theorem~\ref{th:sat-incl-co} is applied to determine ellipsoidal reachable sets. However, due to its generality, the result applies to arbitrary convex sets. 
\end{remark}

The next corollaries follow from Theorem \ref{th:sat-incl-co} by convexity.
\begin{corollary}\label{cor:cor-1}
   \emph{If Assumption \ref{ass:regularity-conds} holds, and $e$ is such that 
\begin{equation}\label{eq:cor1}
e^\top \ABKc^\top P \ABKc e \leq \gamma,
\end{equation}    
   holds for all $\mathscr{J} \subseteq \mathbb{N}^+_m$ with $\gamma \geq 0$, then
   \begin{equation}\label{eq:eqgamma}
      f(e,v)^\top P f(e,v) \leq \gamma, \qquad \forall e \in \R^n.
   \end{equation}}
\end{corollary}
\proof Condition (\ref{eq:cor1}) ensures that every extreme point of the convex hull (\ref{eq:co}) lies in $\mathscr{E}(P,\gamma)$. From Theorem \ref{th:sat-incl-co}, this implies that $f(e,v) \in 
\mathscr{E}(P,\gamma)$, thus inequality (\ref{eq:eqgamma}) holds.
\endproof

\begin{corollary}\label{cor:2norm-bound}
   \emph{If Assumption \ref{ass:regularity-conds} holds, then
      \begin{align}\label{eq:max-ineq}
         \begin{split}
      &f(e,v)^\top P f(e,v) \\
      &\leq \max_{i \in \mathscr{J}} \,\, e^\top \ABKc^\top P \ABKc e,
         \end{split}
   \end{align}
   holds for all $e \in \R^n$ and $i \in \mathscr{J} \subseteq \mathbb{N}^+_m$.}
\end{corollary}

\proof As in Cor. \ref{cor:cor-1}, the claim follows from Theorem \ref{th:sat-incl-co}.\endproof

These findings indicate that the behavior of the nonlinear system $f(e, v)$ is effectively bounded by that of the linear system $\ABKc e$ for all $\mathscr{J} \subseteq \mathbb{N}^+_m$. We can then construct bounds on $\ABKc e$ using linear systems theory which hold for the saturated system.\\

\subsection{Base PRS Construction}

Using the results presented in section \ref{sec:co-bounds-sat} above, we will construct a PRS
for the error $e_k \in \R^n$.

\begin{proposition}\label{prop:base-exp-bound}
   \emph{Let Assumption \ref{ass:regularity-conds} hold, and consider a matrix $P = P^\top \succ 0$ satisfying 
   \begin{equation}\label{eq:contractivity-lmi}
      \begin{split}
      \ABKc^\top P \ABKc \preceq \lambda P,
      \end{split}
   \end{equation}
for all $\J \subseteq \N_m^+$ and contraction rate $\lambda \in [0, 1)$. Then,
   \begin{equation}\label{eq:exp-bound}
      \begin{split}
      \E\left[ e_{k}^\top P e_{k} \right] \leq \frac{1 - \lambda^k}{1 - \lambda} \mathrm{Tr}(PW),
      \end{split}
   \end{equation}
   holds, with $\E\left[ e_0^\top P e_0 \right] = 0$.}
\end{proposition}

\proof If (\ref{eq:contractivity-lmi}) holds, Corollaries \ref{cor:cor-1} and \ref{cor:2norm-bound} yield

\begin{equation*}\label{eq:ellipsoid-bound}
   f(e_k, v_k)^\top P f(e_k, v_k) \leq \lambda e^\top_k P e_k,
\end{equation*}
which gives the relationship below
\begin{align*}
   \begin{split}
      &\E\left[ \left(f(e_k, v_k) + w_k\right)^\top P \left(f(e_k, v_k) + w_k\right) \right] \\
      &\qquad = \E\left[ f(e_k, v_k)^\top P f(e_k, v_k) \right] + \E\left[w_k^\top P w_k \right] \\
      &\qquad = \E\left[ f(e_k, v_k)^\top P f(e_k, v_k) \right] + \E\left[\mathrm{Tr}\left(w_k^\top P w_k\right) \right], \\
      &\qquad = \E\left[ f(e_k, v_k)^\top P f(e_k, v_k) \right] + \mathrm{Tr} \left(P W \right), \\
      &\qquad \leq \lambda\E\left[ e_k^\top P e_k \right] + \mathrm{Tr} \left(P W \right).
   \end{split}
\end{align*}

Given the error dynamics (\ref{eq:error-dynamics}), we have
\begin{equation}\label{eq:recursive-exp-bound}
   \E\left[e_{k+1}^\top P e_{k+1} \right] \leq \lambda \E\left[e_k^\top P e_k \right] + \mathrm{Tr}\left( PW \right).
\end{equation}
The rest of the proof follows by induction. Recall that for $k = 0$, we have $\E\left[ e_0^\top P e_0 \right] = 0$. For $k = 1$, we have 
\begin{align*}
\begin{split}
    \E\left[ e_1^\top P e_1 \right] &\leq \lambda \E\left[e_0^\top P e_0 \right] + \mathrm{Tr}(PW), \\
    &\leq \mathrm{Tr}(PW).
\end{split}
\end{align*}
For $k = t$, we have 
\begin{align*}\label{eq:geom-series}
\begin{split}
    \E\left[ e_t^\top P e_t \right] &\leq \lambda \E\left[e_{t-1}^\top P e_{t-1} \right] + \mathrm{Tr}(PW), \\
    &\leq \sum_{i = 1}^{t} \lambda^{i} \mathrm{Tr}(PW). \\
\end{split}
\end{align*}
The successive upper-bounds on $\E\left[e_{k}^\top P e_{k}\right]$ represent a geometric series with common 
ratio $\lambda$. Therefore, (\ref{eq:exp-bound}) holds for any contraction rate $0 \leq \lambda < 1$. \endproof

\begin{remark}
    While the simplifying i.i.d. assumption \ref{ass:noise-iid} is fairly standard, it can often be unrealistic. The above bound does not hold in the presence of biased or correlated disturbances, which would introduce the additional terms $\E \left[w_k \right]^\top P \E \left[w_k \right]$ and $\E \left[ f \left( e_k, v_k \right)^\top P w_k  \right]$, respectively. As in \cite{fiacchiniProbabilistic2021}, these can be dealt with by assuming that bounds on the means and covariance matrices exist. 
\end{remark}

Using the above result in combination with the Markov bound, we present the following Lemma.

\begin{lemma}\label{lem:markov-reach}
   \emph{The sequence
   \begin{equation}\label{eq:k-step-markov}
      \mathcal{R}^\varepsilon_k(\lambda) = \left\{e \; \middle| \; e^\top P e \leq 
      \frac{1 - \lambda^k}{\varepsilon \left( 1 - \lambda \right)}\mathrm{Tr}\left( PW \right) \right\},
   \end{equation}
   is a Markov-based ellipsoidal PRS associated with the violation level $\varepsilon$.}
\end{lemma}

\proof The claim follows from the application of Markov's inequality on the random variable 
$e_{k}^\top P e_{k}$
\begin{equation*}
   \Pr\left\{ e_{k}^\top P e_{k} \geq r \right\} \leq \frac{\E\left[ e_{k}^\top P e_{k} \right]}{r},
\end{equation*}
and combined with (\ref{eq:exp-bound}) for $r = \frac{1 - \lambda^k}{\varepsilon \left( 1 - \lambda \right)}\mathrm{Tr}\left( PW \right)$, giving,
\begin{align*}
   \Pr\left\{ e_{k}^\top P e_{k} \leq r \right\} &\geq 1 - \frac{1 - \lambda^k}{r\left(1 - \lambda\right)} \mathrm{Tr}(PW), \\
   & = 1 - \varepsilon. \tag*{\QED}
\end{align*}

Computing the Markov-based PRS $\mathcal{R}^\varepsilon_k (\lambda)$ for $k \to \infty$ yields the Markov-based PUB defined
below.

\begin{lemma}\label{lem:markov-bound}
   \emph{The set
   \begin{equation}
      \mathcal{R}^\varepsilon(\lambda) = \left\{e \; \middle | \; e^\top P e \leq 
      \frac{1}{\varepsilon \left( 1 - \lambda \right)}\mathrm{Tr}\left( PW \right) \right\},
   \end{equation}
   is a Markov-based ellipsoidal PUB associated with the violation level $\varepsilon$.}
\end{lemma}

\proof The claim follows from Lemma \ref{lem:markov-reach} and the fact that
\begin{equation*}
   \lim_{k \to \infty} \frac{1 - \lambda^k}{\left(1 - \lambda\right)} \mathrm{Tr}(PW) 
   = \frac{1}{1 - \lambda}\mathrm{Tr}\left( PW \right). \tag*{\QED}
\end{equation*}

These results provide a method to compute well-defined PRS for a saturated system at any time step $k$ while guaranteeing that they will converge to a finite geometry as $k \to \infty$. This property can be leveraged for control applications as it allows us to establish conditions under which these sets remain bounded within the original constraint set $\mathcal{X}$.

It is nevertheless important to note that the bound derived in this section, which aligns with the one presented in \cite{kohlerPredictive2025} for a linearized system, can be overly conservative.  In fact, condition (\ref{eq:contractivity-lmi}) directly implies
\begin{equation}\label{eq:ol-contraction-lmi}
      A^\top P A \preceq \lambda P,
\end{equation}
meaning that $\lambda$ is a contraction rate valid for the open-loop and, equivalently, the saturated closed-loop systems. This bound therefore neglects the improved performance of the system near the origin, where the non-saturated closed-loop dynamics dominate. This conservatism can be especially significant when the saturation bounds are high. To address this limitation, the following section introduces a refined bound, which constitutes the main contribution of this paper.

\subsection{Conditional Tightening of Base PRS}
The Markov-based PRS and PUB presented in Lemmas \ref{lem:markov-reach} and \ref{lem:markov-bound} are unnecessarily conservative,
as they rely on a contraction rate $\lambda$ that does not fully capture the system dynamics. In fact, there exists a region in the state space where system (\ref{eq:error-dynamics}) reduces to the closed-loop system (where the control input $Ke + v$ does not saturate), called the region of linearity. An expression of that region is given below.

\begin{definition}[\cite{tarbouriechStability2011}]
   The region of linearity of system (\ref{eq:error-dynamics}), denoted $\mathcal{R}_L$, is defined as the
   set of all error terms such that $\varphi\left(Ke + v\right) = Ke + v$.
\end{definition}

Considering that the error dynamics (\ref{eq:error-dynamics}) reduce to $e_{k+1} = (A+ BK)e_k + w_k$ for $e_k \in \mathcal{R}_L$, we present the corollary below which follows from Theorem \ref{th:sat-incl-co}, Proposition \ref{prop:base-exp-bound},
and the definition presented above.

\begin{corollary}\label{cor:lambdal}
   \emph{Let Assumptions \ref{ass:noise-iid} and \ref{ass:regularity-conds} 
   hold, and consider the matrix $P = P^\top \succ 0$ which solves the Lyapunov equation 
   (\ref{eq:contractivity-lmi}), with contraction rate $\lambda \in [0, 1)$. Suppose that the gain
   $K \in \mathbb{R}^{n \times m}$ and $\lambda_L \in \mathbb{R}$ are such that
   \begin{equation}\label{eq:find-lambdal}
      \begin{split}
   &      \lambda_L < \lambda,\\
   &      \begin{bmatrix}
            \lambda_L P & (A + BK)^\top \\
            A + BK & P^{-1}
         \end{bmatrix} \succeq 0,
      \end{split}
   \end{equation}
   then,
   \begin{equation}\label{eq:lambdal-bound}
      \Ek{k+1} \leq \lambda_L \Ek{k} + \Tr(PW),
   \end{equation}
   for all $e_k \in \mathcal{R}_L$ and all $k \in \mathbb{N}_0$.}
\end{corollary}

\proof Condition (\ref{eq:find-lambdal}) states that the contraction rate $\lambda_L$, valid in the region $\mathcal{R}_L$, is smaller than $\lambda$. Indeed, by the Schur complement, it is equivalent to
\begin{equation}\label{eq:ABK}
   (A + BK)^\top P (A + BK) \preceq \lambda_L P.
\end{equation}
Thus, when $e_k \in \mathcal{R}_L$, (\ref{eq:recursive-exp-bound}) holds for $\lambda_L$, yielding (\ref{eq:lambdal-bound}). 
\endproof
To reconcile the open-loop contraction \eqref{eq:contractivity-lmi} with the local closed-loop contraction valid in $\mathcal{R}_L$, we impose the following compatibility condition.
\begin{assumption}\label{ass:stable-cl}
   Matrices $K$ and $P$ and parameters $\lambda$ and $\lambda_L$ are such that conditions (\ref{eq:contractivity-lmi}) and (\ref{eq:find-lambdal}) hold.
\end{assumption}
\begin{remark} 
Such an assumption entails no loss of generality, as it simply means that the open-loop dynamics converge more slowly than the closed-loop dynamics. In general, a gain $K$ can be explicitly designed (by solving a suitable optimization problem) to ensure this assumption holds for a common Lyapunov function $P$.
\end{remark}

Corollary \ref{cor:lambdal} provides a bound on the evolution of $\mathbb{E}\left[e_k^\top P e_k\right]$ that holds whenever $e_k$ lies in the region of linearity, while (\ref{eq:contractivity-lmi}) gives a more conservative bound that holds over the entire state space. Taken together, these two results allow us to exploit distributional information about the error term to reduce the conservatism of sets (\ref{eq:k-step-markov}). Specifically, we can more accurately estimate $\E\left[e_{k}^\top P e_{k}\right]$ when considering that the system frequently operates in $\mathcal{R}_L$ where $\lambda_L$ governs the dynamics. To formalize this within our ellipsoidal framework, we state the following theorem.

\begin{theorem}[\cite{tarbouriechStability2011}]\label{th:ellipsoidal-Rl}
   \emph{The maximal ellipsoid with shape matrix $P$ contained in the region of linearity $\mathscr{E}(P, r_L) \subseteq \mathcal{R}_L$ has a scaling parameter $r_L$ defined by
   \begin{equation}
      r_L = \min_{1 \leq i \leq m} \frac{\left(1 - v_i \right)^2}{K_{i} P^{-1} K_{i}^\top}.
   \end{equation}}
\end{theorem}
We will utilize this ellipsoid to determine whether the error terms are in the region of linearity $\mathcal{R}_L$, i.e., if the error is such that $e^\top P e \leq r_L$, then $\varphi\left(Ke + v\right) = Ke + v$.

\begin{proposition}\label{prop:balanced-bound}
   \emph{Let Assumptions \ref{ass:noise-iid}, \ref{ass:regularity-conds} and \ref{ass:stable-cl} 
   hold. 
   Then,
   \begin{equation}\label{eq:exp-bound-tight}
      \begin{split}
         \Ek{k+1} \leq &\lambda_L \E\left[e_k^\top P e_k \right] + \frac{\lambda - \lambda_L}{r_L} \E^2\left[e_k^\top P e_k \right]
         \\ &+ \Tr(PW).
      \end{split}
   \end{equation}
   holds for all $k \in \mathbb{N}_0$.}
\end{proposition}

\proof 
Since we have distinct contraction rates $\lambda_L$ and $\lambda$ for 
$e_k \in \mathcal{R}_L$ and $e_k \notin \mathcal{R}_L$ respectively, 
we can partition the state space $\R^n$ into $\mathcal{S}_1 = \set{e \in \R^n : e^\top P e \leq r_L}$ and $\mathcal{S}_2 = \set{e \in \R^n : e^\top P e > r_L}$, where, by construction, $S_1 \subseteq \mathcal{R}_L$. Denote $q_{k+1} = e_{k+1}^\top P e_{k+1}$. By the law of total expectation
\begin{align*}
   \begin{split}
   \E\left[ q_{k+1} \right] &= \sum_{i = 1}^2 \E\left[q_{k+1}\mid 
   e_k \in \mathcal{S}_i\right] \Pr\left\{e_k \in \mathcal{S}_i\right\}.
   \end{split}
\end{align*}
Given that $\Pr \left\{e_k \in S_2 \right\} = 1 - \Pr \left\{e_k \in \mathcal{S}_1 \right\}$, we can write
\begin{equation*}
   \begin{split}
      \Eq{k+1} = & \big(\E\left[q_{k+1} \mid e_k \in \mathcal{S}_1\right] - \E\left[q_{k+1} \mid e_k \in  \mathcal{S}_2\right]\big) \\
      & \cdot \Pr \left\{e_k \in  \mathcal{S}_1 \right\} + \E \left[ q_{k+1} \mid e_k \in \mathcal{S}_2 \right].
   \end{split}
\end{equation*}
Combining the bounds (\ref{eq:recursive-exp-bound}) and (\ref{eq:lambdal-bound}) applied for $\lambda$ and $\lambda_L$ with the Markov bound on $\Pr \left\{e_k \in \mathcal{S}_1 \right\}$ one obtains
\begin{subequations}\label{eq:exp-bounds-tight}
   \begin{align}
      \begin{split}\label{eq:upper-bound-form}
         \Eq{k+1}&\leq \lambda\Eq{k} + \left( \lambda_L - \lambda \right) \Eq{k} \left( 1 - \frac{\Eq{k}}{r_L} \right)  \\
         &\qquad+ \Tr(PW),
         \end{split} \\
         &\leq \lambda_L \Eq{k} + \frac{\lambda - \lambda_L}{r_L} \E^2\left[ q_k \right] + \Tr(PW).
   \end{align} 
\end{subequations}
Therefore, (\ref{eq:exp-bound-tight}) holds for any pair of contraction rates $\left(\lambda_L, \lambda\right)$ such that $0 \leq \lambda_L < \lambda < 1$.
\endproof

Proposition \ref{prop:balanced-bound} states that both contraction rates $\lambda_L$ and $\lambda$ contribute to the bound on $\Ek{k}$. The interaction between these rates suggests the existence of an effective contraction rate that appropriately averages their contributions. Such an effective rate would enable us to construct a more interpretable bound on the expectation, similar to (\ref{eq:exp-bound}). This motivates our subsequent analysis to identify and characterize this average behavior. Note that in what follows, we identify a contraction rate valid for all $k \in \mathbb{N}_0$, but the analysis can be easily restricted to some finite value of $k$.

\begin{proposition}\label{prop:tighter-bound-on-exp}
   \emph{Let Assumptions \ref{ass:noise-iid}, \ref{ass:regularity-conds} and \ref{ass:stable-cl} hold, and consider bound 
   (\ref{eq:exp-bound-tight}).
   If 
   \begin{equation}\label{eq:rl}
    \frac{1}{1 - \lambda} \Tr(PW) < r_L,
   \end{equation}
   then
   \begin{equation}\label{eq:exp-bound-bar-cf}
      \Ek{k} \leq \frac{1 - \left(\bar\lambda^\ast\right)^k}{1 - \bar\lambda} \Tr(PW),
   \end{equation}
   holds for 
   $\bar\lambda^\ast \in \left[ \lambda_L,\; \lambda \right]$ such that  
   \begin{equation}\label{eq:lambdastar}
         \frac{1}{1 - \bar\lambda^\ast} \Tr(PW) = \frac{\bar\lambda^\ast - \lambda_L} {\lambda - \lambda_L} r_L.
   \end{equation}}
\end{proposition}

\proof 
Since $\lambda_L - \lambda \leq 0$, inequality (\ref{eq:upper-bound-form}) implies
\begin{equation*}
   \Ek{k+1} \leq \lambda \Ek{k} + \Tr(PW).
\end{equation*}
Iterating this recursion yields
\begin{equation*}
   \begin{split}
     \Ek{k} &\leq \frac{1 - \lambda^k}{1 - \lambda} \Tr(PW) \leq \frac{1}{1 - \lambda} \Tr(PW),
   \end{split}
\end{equation*}
for all $k \in \mathbb{N}_0$ and every $e_k \in \R^n$. From condition (\ref{eq:rl}) there exists a value $\mu \in \left[\lambda_L,\; \lambda \right]$ such that
\begin{equation}\label{eq:exp/rl-param}
    \frac{\Ek{k}}{r_L} \leq \frac{\mu - \lambda_L}{\lambda - \lambda_L} \leq 1.
\end{equation}
Combining this parametrized bound with (\ref{eq:exp-bound-tight}) gives
\begin{equation}\label{eq:exp-bound-bar}
   \Ek{k+1} \leq \mu \Ek{k} + \Tr(PW),
\end{equation}
which is exactly (\ref{eq:recursive-exp-bound}) for parameter $\mu$. Iterating yields
\begin{equation}\label{eq:bound-bar}
      \Ek{k} \leq \frac{1 - \mu^k}{1 - \mu} \Tr(PW) \leq \frac{1}{1 - \mu} \Tr(PW),
\end{equation}
provided (\ref{eq:exp/rl-param}) holds. 
We now seek the values of $\mu$ for which (\ref{eq:exp/rl-param}) holds for all $k \in \N_0$ and thus ensures (\ref{eq:bound-bar}). Define
\begin{equation}
   g(\mu) \coloneqq \frac{\mu - \lambda_L}{\lambda - \lambda_L} r_L - \frac{1}{1 - \mu} \Tr(PW),
\end{equation}
a continuous and strictly concave function on the interval $\left[ \lambda_L, \lambda \right]$. From (\ref{eq:rl}), we have $g(\lambda) > 0$. At $\mu = \bar\lambda^\ast$, defined by (\ref{eq:lambdastar}), we have $g(\bar{\lambda}^\ast) = 0$. By strict concavity, this implies $g(\mu) > 0$ for all $\mu \in (\bar{\lambda}^\ast, \lambda]$. Equivalently,
\begin{equation}\label{eq:mu}
\frac{1}{1 - \mu} \Tr(PW) < \frac{\mu - \lambda_L}{\lambda - \lambda_L} r_L, \quad \forall \mu \in (\bar{\lambda}^\ast, \lambda].
\end{equation}
Fix $\mu \in (\bar{\lambda}^\ast, \lambda]$. Three cases are possible for $\Ek{k}$:
\begin{enumerate}
   \item $\displaystyle{\Ek{k} \leq \frac{1}{1 - \mu} \Tr(PW),}$ \label{case1}
   \item $\displaystyle{\frac{1}{1 - \mu} \Tr(PW) < \Ek{k} \leq \frac{\mu - \lambda_L}{\lambda - \lambda_L} r_L,}$ \label{case2}
   \item $\displaystyle{\frac{\mu - \lambda_L}{\lambda - \lambda_L} r_L < \Ek{k}}$. \label{case3}\\
\end{enumerate}

Case \ref{case2}) cannot occur, as it contradicts the fact that (\ref{eq:exp/rl-param}) implies (\ref{eq:bound-bar}).

Case \ref{case3}) cannot occur: by concavity of $g$, if $\Ek{k}$ exceeded $\frac{\mu - \lambda_L}{\lambda - \lambda_L} r_L$, there would exist a $\nu \in (\mu,\lambda]$ for which  
\begin{equation*}
   \frac{1}{1-\nu}\Tr(PW) < \Ek{k} = \frac{\nu - \lambda_L}{\lambda - \lambda_L}r_L,    
\end{equation*}
where the left-hand side follows from (\ref{eq:mu}). This falls into Case \ref{case2}, a contradiction.  


Therefore, only case \ref{case1}) is possible. Analogous considerations also hold for $\mu = \bar{\lambda}^\ast$.

This means that any $\mu \in [\bar\lambda^\ast,\lambda]$ is such that (\ref{eq:bound-bar}) holds for all $k \in \N_0$. The smallest such $\mu = \bar\lambda^\ast$ yields the tightest valid upper bound:  
\[
   \Ek{k} \leq \frac{1 - (\bar\lambda^\ast)^k}{1 - \bar\lambda^\ast} \Tr(PW).
\]  
This is precisely (\ref{eq:exp-bound-bar-cf}), which concludes the proof. \endproof

The result in Proposition \ref{prop:tighter-bound-on-exp} indicates that, under certain conditions,
mainly when $\mathscr{E}(P, r_L)$ is sufficiently large, there exists a $\bar\lambda^\ast < \lambda$ that allows us to construct a deterministic upper bound $\Ek{k}$ that is tighter than the one given in (\ref{eq:exp-bound}) and larger than that in (\ref{eq:lambdal-bound}). However, when $r_L$ is smaller than or equal to the critical value $\frac{1}{1 - \lambda}\Tr(PW)$, we can no longer ensure that $\Ek{k} \leq r_L$ for all $k \in \mathbb{N}_0$. Then, the Markov inequality used to construct these bounds
\begin{equation}
   \Pr \left\{e_k^\top P e_k \leq r_L \right\} \geq 1 - \frac{\Ek{k}}{r_L},
\end{equation}
no longer provides any meaningful information on the probability that the error term is in the region of linearity. In such a case, the only admissible bound on the expectation of the error
is that provided in (\ref{eq:exp-bound}). More specifically, we have
\begin{equation}\label{eq:hatlambda}
    \hat\lambda = \begin{cases}
        \bar\lambda^\ast &\qquad \text{if } \frac{1}{1 - \lambda} \Tr(PW) < r_L, \\
        \lambda &\qquad \text{if } \frac{1}{1 - \lambda} \Tr(PW) \geq r_L,
    \end{cases}
\end{equation}
where $\bar{\lambda}^\ast < \lambda$ is given by (\ref{eq:lambdastar}).

\begin{lemma}\label{lem:markov-reach-tight}
   \emph{The sequence
   \begin{equation}
      \mathcal{R}^\varepsilon_k(\hat\lambda) = \left\{e \; \middle | \; e^\top P e \leq 
      \frac{1 - \hat\lambda^k}
      {\varepsilon \left( 1 - \hat\lambda \right)}\mathrm{Tr}\left( PW \right) \right\},
   \end{equation}
   is a sequence of Markov-based ellipsoidal PRS associated with the violation level $\varepsilon$,
   with $\hat\lambda$ as in (\ref{eq:hatlambda}).}
\end{lemma}

\proof The claim follows from Proposition \ref{prop:tighter-bound-on-exp} and by applying Markov's inequality to the random variable $e_{k+1}^\top P e_{k}$.
\endproof

Computing the Markov-based PRS $\mathcal{R}^\varepsilon_k(\hat\lambda)$ for $k \to \infty$ yields the Markov-based PUB defined below.

\begin{lemma}\label{lem:markov-bound-tight}
   \emph{The set
   \begin{equation}
      \mathcal{R}^\varepsilon(\hat{\lambda}) = \left\{e \; \middle | \; e^\top P e \leq 
      \frac{1}{\varepsilon \left( 1 - \hat\lambda \right)}\mathrm{Tr}\left( PW \right) \right\},
   \end{equation}
   is a Markov-based ellipsoidal PUB associated with the violation level $\varepsilon$, with $\hat\lambda$ as in (\ref{eq:hatlambda}).
   }
\end{lemma}

The sets introduced above are analogous to those defined earlier in Lemmas \ref{lem:markov-reach} and \ref{lem:markov-bound}, but with a key difference in their construction. They primarily employ the bound in (\ref{eq:bound-bar}) while maintaining (\ref{eq:exp-bound}) as a conservative fallback bound when the conditions are unsuitable for computing the effective contraction rate $\bar\lambda^\ast$. The next section will illustrate this approach through a numerical example.

\section{Numerical Example}
We demonstrate the applicability of the proposed method for estimating reachable sets on a simple system. The code is available online\footnote{https://github.com/CarloKaram/PRS-Sat-Sys}.

Consider system (\ref{eq:lti-sat-system}) with
\begin{equation*}
  A = \begin{bmatrix} 0.89 & 0.1 \\ 0.1 & 0.89 \end{bmatrix}, \quad B = \begin{bmatrix} 0 \\ 1 \end{bmatrix}, \quad W = \begin{bmatrix} 1 & 0 \\ 0 & 1 \end{bmatrix},
\end{equation*} 
and saturation bounds ${u_\text{min}, u_\text{max}} = \pm 10$. The noise is generated from a normal
distribution $w_k \sim \mathcal{N}\left(0, W\right)$. For simplicity, we consider $e_k = x_k$ and $v_k = 0$ for all $k \in \N_0$. Considering a dynamic $v_k \neq 0$ (i.e., in a closed-loop SMPC simulation) introduces additional complexities and requires a more involved treatment. As such, we defer that to a future work and focus here on demonstrating the benefit of our method for a fixed $v_k = 0$ for all $k \in \N_0$.  

For a gain $K = \begin{bmatrix} -0.282 & -0.8415\end{bmatrix}$, Assumption \ref{ass:stable-cl} is indeed met. In fact, we solve an optimization problem satisfying (\ref{eq:contractivity-lmi}) while minimizing $\lambda$ via bisection, and obtain an open-loop contraction rate $\lambda = 0.9802$ and a Lyapunov function $P = \begin{bmatrix} 3.54 & 0.67 \\ 0.67 & 3.51 \end{bmatrix}$. Using the obtained value for $P$, we then solve another optimization problem which satisfies (\ref{eq:find-lambdal}) while minimizing $\lambda_L$, and get an unsaturated closed-loop contraction rate $\lambda_L =0.7684$. Finally, we construct the bound (\ref{eq:exp-bound-bar-cf}) using the computed values and design the corresponding PRS and PUB at violation level $\varepsilon = 0.2$.

The maximal inscribed ellipsoid $\mathscr{E}(P, r_L) \subseteq \mathcal{R}_L$ resulting from these parameters has a scaling parameter $r_L = 485.47$,
which yields a value $\hat\lambda = \bar\lambda^\ast = 0.7826$. Figure \ref{fig:exp-bound-comparison} shows that the bound constructed using $\bar\lambda^\ast$ is very  close to the one based on $\lambda_L$, both resulting in rather tight upper estimations of the true expectation, calculated from $10^3$ noise samples. In contrast, the bound derived from $\lambda$ is overly conservative and grows excessively large, failing to converge within $k = 100$. Figure \ref{fig:pub-comparison} shows that the PUB
$\mathcal{R}^{0.2}\left(\bar\lambda^\ast\right)$ is significantly smaller than $\mathcal{R}^{0.2}(\lambda)$. In fact, the reduction in scaling from $r_\lambda \approx 1777$ to $r_{\bar\lambda^\ast} \approx 162$, corresponds to a reduction in PUB conservatism (area) of approximately $91\%$ with no compromise on probabilistic containment guarantees.

\addtolength{\textheight}{-0.243cm}

Finally, Figure \ref{fig:lambdabar-evolution} illustrates how the effective contraction rate $\bar\lambda^\ast$ approaches $\lambda_L$ as $r_L$ increases. The range of $r_L$ shown corresponds to symmetric saturation bounds $\left|u_\text{min}\right| = \left|u_\text{max}\right| > 8.5$, below which condition (\ref{eq:rl}) is not met and (\ref{eq:exp-bound}) becomes the only admissible bound for all $k$. Clearly, a larger region of linearity reduces the likelihood of input saturation, allowing for more aggressive control in an essentially linear domain. This result emphasizes the significance of accounting for the region of linearity when analyzing reachability in the presence of input saturation, and further highlights the effectiveness of this approach in reducing the conservatism of existing PRS and PUB construction methods.

\begin{figure*}[hb]
   \vspace{0.3cm}
   \centering
   \begin{subfigure}{0.66\textwidth}
   \captionsetup{skip=0pt}
      \centering
      \begin{tikzpicture}
         \begin{groupplot}[group style={group size=3 by 1, horizontal sep=0.8cm},small]
         \nextgroupplot[
            xlabel={$k$},
            ylabel={$\Ek{k}$},
            legend style={font=\tiny},
            legend cell align=left,
            legend pos=north west,
         ]
            \addplot[
               thin, 
               solid, 
               red, 
               mark=x,
               mark repeat=5,
               mark phase=1,
               ] table [x=k, y=e, col sep=comma] {data.csv};
            \addplot[black, thick, solid] table [x=k, y=l, col sep=comma] {data.csv};
            \addplot[green!60!black, thick, solid,
            ] table [x=k, y=ll, col sep=comma] {data.csv};
            \addplot[blue, thick, solid] table [x=k, y=lb, col sep=comma] {data.csv};
            \coordinate (nw) at (15, 40);
            \coordinate (sw) at (15, 10);
            \coordinate (ne) at (101, 40);
            \coordinate (se) at (101, 10);
            \draw[thin, solid] (15,10) rectangle (101,40);
            \legend{True Mean, $\lambda$ bound, $\lambda_L$ bound, $\bar\lambda^\ast$ bound};
         \nextgroupplot[
            restrict x to domain=15:100,
            ytick distance = 4]
            \addplot[
               thin, 
               solid, 
               red, 
               mark=x,
               mark repeat=5,
               mark phase=1,
               ] table [x=k, y=e, col sep=comma] {data.csv};
               \addplot[green!60!black, thick, solid] table [x=k, y=ll, col sep=comma] {data.csv};
               \addplot[blue, thick, solid] table [x=k, y=lb, col sep=comma] {data.csv};
         \end{groupplot}
            \draw[dotted] (nw) -- (group c2r1.north west);
            \draw[dotted] (sw) -- (group c2r1.south west);
            \path[name path=lineA] (se) -- (group c2r1.south east);
            \path[name path=lineB] (ne) -- (group c2r1.north east);
            \path[name path=lineC] (group c2r1.north west) -- (group c2r1.south west);
            \draw [dotted,name intersections={of=lineA and lineC, by=pA}] (se) -- (pA);
            \draw [dotted,name intersections={of=lineB and lineC, by=pB}] (ne) -- (pB);
            \draw [dotted,gray!90,opacity=0.7] (pA) -- (group c2r1.south east);
            \draw [dotted,gray!90,opacity=0.7] (pB) -- (group c2r1.north east);
      \end{tikzpicture}
      \caption{}\label{fig:exp-bound-comparison}
   \end{subfigure}
   \begin{subfigure}{0.33\textwidth}
      \captionsetup{skip=0pt}
      \centering
      \begin{tikzpicture}
         \begin{axis}[
              small,
              xlabel={$x_1$},
              ylabel={$x_2$},
              legend style={font=\tiny},
              legend cell align=left,
              legend pos=north west,
              ylabel shift=-8pt,
          ]
          \addplot[black, thick, solid] table [x=x, y=y, col sep=comma] {lell.csv};
          \addplot[blue, thick, solid] table [x=x, y=y, col sep=comma] {lbell.csv};
          \addplot[
          scatter,
          only marks,
         point meta=explicit symbolic,
          scatter/classes={
              a={mark=*, red} 
          },
          mark options = {scale=0.5},
          red,
          ] table[x=x, y=y, col sep=comma] {states.csv};
          \legend{$\lambda$-based PUB, $\bar\lambda^\ast$-based PUB, True States}
          \end{axis}
      \end{tikzpicture}
      \caption{}\label{fig:pub-comparison}
  \end{subfigure}
  \caption{Conditional bound tightening applied to the system -- (a) comparison of the bounds on $\Ek{k}$ derived from $\lambda$, $\lambda_L$ and $\bar\lambda^\ast$ with the sampled expectation; (b) comparison of $\lambda$ and $\bar\lambda^\ast$-based PUBs and their error sample containment.}
  \vspace{-0.4cm}
\end{figure*}

\begin{figure}[ht]
   \centering
   \begin{tikzpicture}[scale=0.85]
       \begin{axis}[
           xlabel={$r_L$},
           ylabel={$\bar\lambda$},
           legend style={font=\footnotesize},
           legend cell align=left,
           legend pos=north east,
           width=8cm,
           ytick={0.770, 0.775, 0.780, 0.785, 0.790},
           yticklabels={0.770, 0.775, 0.780, 0.785, 0.790},
       ]
       \addplot[dashed, thick, green!60!black] table [x=rl, y=ll, col sep=comma] {convergence.csv};
       \addplot[thick, blue, solid, mark=x] table [x=rl, y=lb, col sep=comma] {convergence.csv};
       \legend{$\lambda_L$, $\bar\lambda^\ast$}
       \end{axis}
   \end{tikzpicture}
   \caption{Evolution of $\bar\lambda^\ast$ with respect to the scaling factor $r_L$.}\label{fig:lambdabar-evolution}
   \vspace{-0.5cm}
\end{figure}


\section{Conclusion}
We presented a constructive method to design probabilistic reachable sets and ultimate bounds in the presence of hard input constraints on a linear system subject to potentially unbounded additive disturbances. By adopting a saturated dynamics formulation, we derive contraction bounds allowing the construction of tight PRS and PUB. Numerical results validated the effectiveness of this novel approach, highlighting its ability to integrate varying constraints. Future work will focus on integrating this technique into SMPC schemes, and extending it to the construction of semialgebraic sets using the Moment-SOS Hierarchy, which incorporates higher-order moments of the noise into the reachability analysis.

\printbibliography

\end{document}